\documentclass[12pt]{amsart}
\usepackage{amsfonts}
\usepackage{graphicx}

  \theoremstyle{plain}
\newtheorem{theorem}{Theorem}[section]
\newtheorem{lemma}[theorem]{Lemma}

\theoremstyle{definition}
\newtheorem{definition}[theorem]{Definition}
\theoremstyle{remark}

\newtheorem{remark}[theorem]{Remark}
\numberwithin{equation}{section}

  \renewcommand\epsilon\varepsilon
  \renewcommand\phi\varphi

\begin{document}
\author{Norayr Matevosyan}
\address{Norayr Matevosyan\\
Johann Radon Institute for Computational and Applied
         Mathematics (RICAM)\\
Austrian Academy of Sciences\\
Altenbergerstraße 69 \\
A-4040 Linz, Austria }
\email{norayr.matevosyan@oeaw.ac.at}
\title[Tangential touch between Free and Fixed boundaries]
{Tangential touch between Free and Fixed boundaries IN A
 PROBLEM FROM SUPERCONDUCTIVITY}
\maketitle

\begin{abstract}
In this paper we study regularity properties of the free
boundary problem
\begin{equation*}
\Delta u=\chi_{\{|\nabla u|\neq0\}} \ \text{in}\ B^+_1,
   \quad u=0\ \text{on}\ B_1\cap\{x_1=0\},
\end{equation*}
where $B^+_1=\{|x|<1,x_1>0\}$ and $B_1=\{|x|<1\}$. If the origin
is a free boundary point, then we show that the free boundary
touches the fixed boundary $\{ x_1=0\}$ tangentially.
\end{abstract}


\section{Introduction}

The aim of this paper is to analyze the regularity of solutions
and the behavior of the free boundary near the fixed one  for a
certain type of free boundary problem. Mathematically the problem
is formulated as follows. Suppose we are given a function $u$ such
that
\begin{equation}
\left\{
\begin{array}{ll}
\Delta u=\chi_{\{|\nabla u|\neq0\}} & \text{in}\ B^+_1,
            \text{in the sense of distributions}\\
u=0 & \text{on}\ \Pi\cap B_1,
\end{array}
\right.  \label{pr}
\end{equation}
where $B^+_1=\{|x|<1,x_1>0\}$, $B_1=\{|x|<1\}$ and
$\Pi=\{x_1=0\}$. Let us denote $\Omega=\Omega(u)=\{|\nabla
u|\neq0\}$, $\Lambda=\Lambda(u)=\{|\nabla u|=0\}$,
$\Gamma(u)=\{x:|\nabla u(x)|=0\}\cap\partial\Omega$ the free
boundary and $\Gamma^*(u)=\Gamma(u)\cap\Pi$ is the set of contact
points (see Figure1).

Note that from classical elliptic regularity theory we have that
the solutions of the problem \eqref{pr} are in the space $C^{1,
\alpha }$, for some $0< \alpha < 1 $.

The main point of interest in this paper is to investigate the
question, ``How do the free and fixed boundaries meet?''

%
\begin{figure}[tbp]
\begin{picture}(300,150)(0,0)
\put(10,0){\includegraphics[height=5cm]{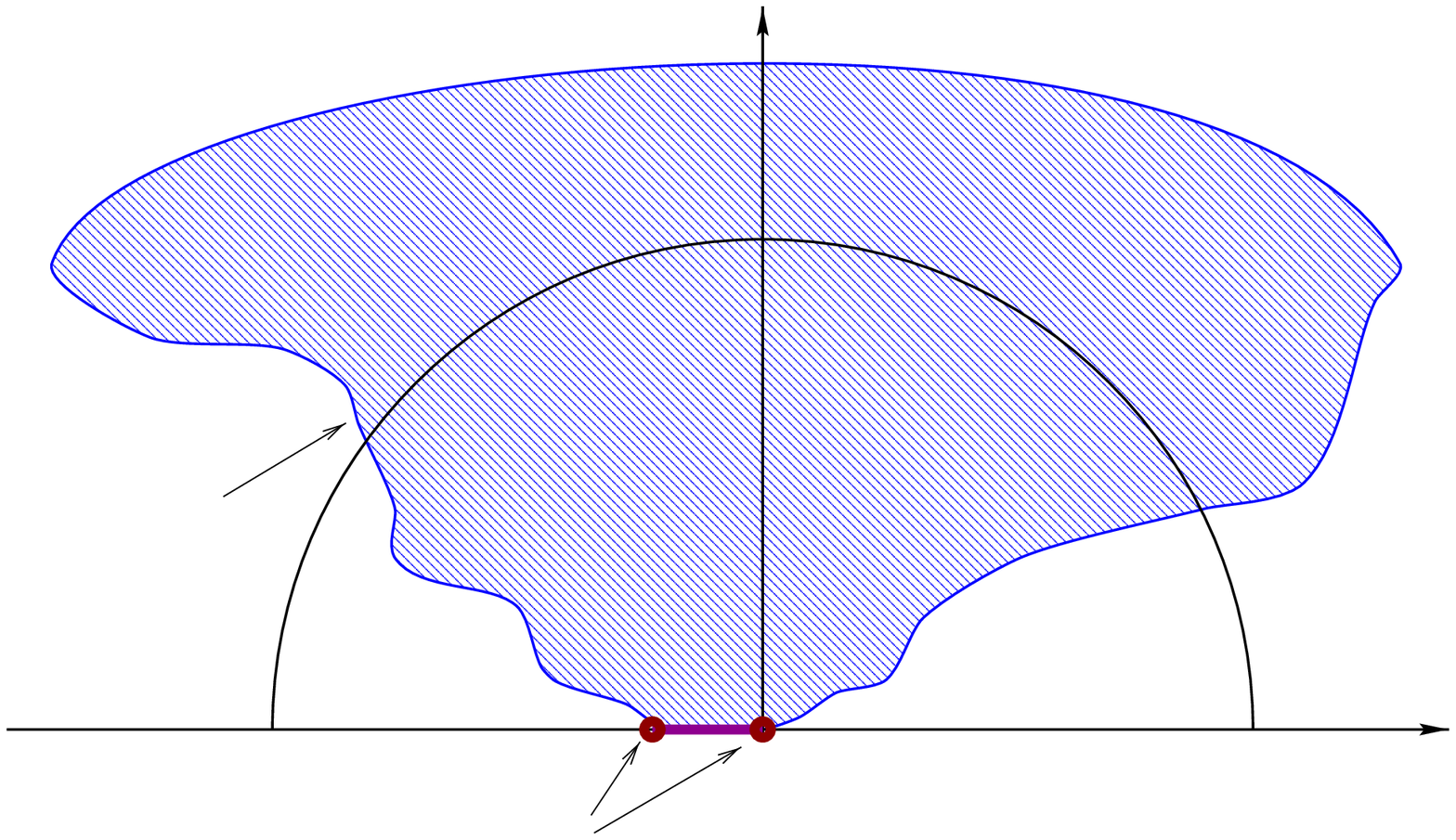}}
\put(160,50){\Large$\Omega$} \put(81,50){$\Delta u=\chi_\Omega$}
\put(182,29){$\Lambda(u)=\{|\nabla u|=0\}$} \put(245,0){$\Pi$}
\put(95,-8){$\Gamma^*$} \put(35,50){$\Gamma$}
\end{picture}
\caption{}
\end{figure}
%
%
\textbf{Notations.} We will use the following notations:\\
\begin{equation*}
\begin{array}{ll}
\mathbf{R}_+^n & \qquad \{x\in\mathbf{R}^n:x_1>0\}, \\
\mathbf{R}_-^n & \qquad \{x\in\mathbf{R}^n:x_1<0\}, \\
B(z,r) & \qquad \{x\in\mathbf{R}^n:|x-z|<r\}, \\
B^+(z,r) & \qquad \{x\in\mathbf{R}^n_+\cap B(z,r)\}, \\
B^-(z,r) & \qquad \{x\in\mathbf{R}^n_-\cap B(z,r)\}, \\
B_r & \qquad B(0,r), \\
B^+_r & \qquad B^+(0,r), \\
B^-_r & \qquad B^-(0,r), \\
\Pi,\ \Pi(z,r),\ \Pi_r & \qquad \{x\in\mathbf{R}^n:x_1=0\},
 \quad \Pi\cap B(z,r),\quad \Pi(0,r), \\
\|\cdot\|_\infty & \qquad \text{supremum norm}, \\
e_1,\ldots,e_n&\qquad \text{standard basis in}\ \mathbf{R}^n,\\
\nu,\ e & \qquad \text{arbitrary unit vectors}, \\
D_\nu,\ D_{\nu e} & \qquad \text{first and second directional
                            derivatives},\\
v_+,\ v_- & \qquad \max(v,0),\ \max(-v,0), \\
\chi_D&\qquad \text{the characteristic function of the set}\ D,\\
\partial D & \qquad \text{the boundary of the set}\ D, \\
\Omega=\Omega(u) & \qquad \{|\nabla u|\neq0\}, \\
\Lambda=\Lambda(u) & \qquad \{|\nabla u|=0\}, \\
\Gamma=\Gamma(u)&\qquad\{x:|\nabla u(x)|=0\}\cap\partial\Omega,\\
\Gamma^*(u) & \qquad \Gamma(u)\cap\Pi.%
\end{array}
\end{equation*}

Free boundary problems, where $\Delta u=\chi_{\{|\nabla
u|\neq0\}}$, appear for instance in connection with
super-conductivity (see [6], [10]). In [8],  [9]  the authors
investigated the problem for the ``interior case,'' i.e. when the
problem is considered in the full ball and there is no fixed
boundary.

A similar problem but with a restriction $u=0$ on $\{|\nabla
u|=0\}$ has been considered earlier by H.Shahgholian and
N.N.Uraltseva in [14]. There the authors have used and developed
further a technique, which mainly uses global analysis as in [5],
[7] and allows one to gain stronger results in problems with no
sign assumption on the solutions.

The elimination of the condition $u=0$ on $\{|\nabla u|=0\}$
generates a number of difficulties in the application of the
technique in [7], [14]. One practical difference is that we no
longer have $u$ vanishing on the free boundary, and this appears in
the technical parts of the proofs. A simple example is that when
scaling we have to take into account the value of the function at
the free boundary points. Also, some of the most crucial tools for
the application of the methods required new proofs. One example is
Lemma \ref{lemma-non-deg}. Another example, which is probably the
most important part of this paper, is the second part of the proof
of Theorem B. Here new geometrical ideas has to be employed, and
these ideas are illustrated in the proof of Lemma \ref{s>0}.

\begin{definition}(Local Solutions)
A function $u$ belongs to the class $P^+_r(M)$, if $u$
satisfies:

\begin{enumerate}
\item $\Delta u=\chi_{\{|\nabla u|\neq0\}}$ in $\ B_r^+$,
      in the sense of distributions,

\item $u=0$ on $\ \Pi_r$,

\item $\|u\|_{\infty,B^+_r}\leq M$.
\end{enumerate}
\end{definition}
Observe that $P^+_r(M)$ is invariant under rotations of
coordinate system, that leave $e_1$ unchanged.
\begin{definition}(Global Solutions)
A function $u$ belongs to the class $P^+_\infty(M)$ if $u$
satisfies:

\begin{enumerate}
\item $\Delta u=\chi_{\{|\nabla u|\neq0\}}$ in $\mathbf{R}_+^n$
      in the sense of distributions,

\item $u=0$ on $\Pi$,

\item $|u(x)|\leq M(|x|+1)^2$.
\end{enumerate}
\end{definition}

We will also need the definition of solution in the whole ball.
\begin{definition}
A function $u$ belongs to the class $P_r(M)$, if $u$ satisfies:

\begin{enumerate}
\item $\Delta u=\chi_{\{|\nabla u|\neq0\}}$ in $\ B_r$,
      in the sense of distributions,

\item $\|u\|_{\infty,B_r}\leq M$.
\end{enumerate}
\end{definition}

 Let us introduce the following notations:
\begin{eqnarray*}
P^+_r(0,M):= & \{u\in P^+_r(M): 0\in\Gamma\},\\
P^+_\infty(0,M):= & \{u\in P^+_\infty(M): 0\in\Gamma\},\\
P_r(0,M):= & \{u\in P_r(M): 0\in\Gamma\}.
\end{eqnarray*}
%
%
%
%
 In the first section we prove $C^{1,1}$-regularity of the
solutions up to $B_{1/2}\cap\Pi$ (see Theorem A) . Then we
classify global solutions in $\mathbf{R}^n_+:=\{x_1>0\}$. Here we
encounter some surprises, in contrast to the problem studied in
[14]. We show that global solutions are either polynomials
(depending on two variables) or one dimensional. The latter will
itself give three different types of solutions (see Theorem B).
Finally we prove our main result, which asserts that the free
boundary touches the fixed one tangentially. \vspace{3mm}

\noindent \textbf{Theorem A}. \textsl{If $u\in P^+_1(0,M)$,
then there is a constant $C=C(n)$ such that
\begin{equation}
\sup_{B^+_{1/2}}{|D_{ij}u|}\leq CM.
\end{equation}
}
\begin{figure}[tbp]
\begin{picture}(400,150)(0,0)
\put(10,-3){\includegraphics[width=12cm]{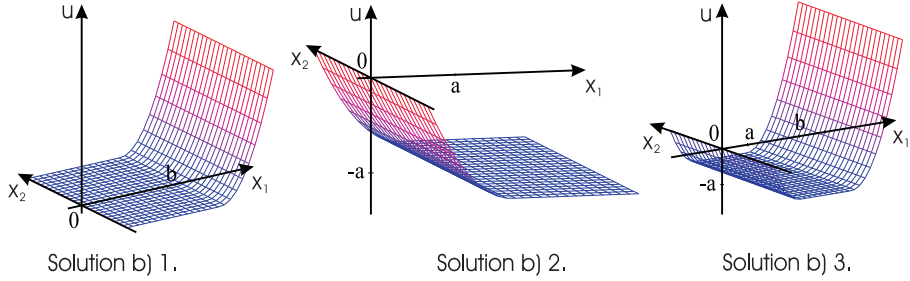}}
\end{picture}
\caption{This figure illustrates three examples of global
solutions for the problem considered in $\mathbf{R}^2$, which we
get under the conditions of Theorem B b), that is for the case
when $\overline{\Omega}\neq\mathbf{R}^2_+$}
\end{figure}

\noindent \textbf{Theorem B}. \textsl{Let $u\in~P_\infty^+(M)$.
Then, in some rotated system of coordinates which leaves $e_1 $
unchanged, the following holds }

\begin{itemize}

\item[a)] \textsl{If $\overline{\Omega}=\overline{\mathbf{R}^n_+}$,
then
\begin{equation*}
u(x)=\frac{x_1^2}2+ax_1x_2+\alpha x_1,\quad \text{with}\
                                       a,\alpha\in\mathbf{R}.
\end{equation*}
}
\item[b)] \textsl{If $\overline{\Omega}\neq\mathbf{R}^n_+$,
  then $u$ depends only on $x_1$ and has one of the following
  representations: }
\begin{enumerate}

\item \textsl{$u(x)=\frac{(x_1-b)_+^2}2$, for $b>0$; }

\item \textsl{$u(x)=\frac{(x_1-a)_-^2-a^2}2$, for $a>0$; }

\item \textsl{$u(x)=\frac{(x_1-a)_-^2+(x_1-b)_+^2-a^2}2$,
      for some $0<a<b$. }
\end{enumerate}
\end{itemize}

\vspace{3mm}

\noindent \textbf{Theorem C}. \textsl{\ There exists
$r_0=r_0(n,M)>0$ and a modulus of continuity
$\sigma(\sigma(0^+)=0)$ such that if $u\in P^+_1(0,M)$,
then}
%
%
\textsl{\
\begin{equation}
\partial\Omega\cap B_{r_0}\subset\{x:x_1\leq\sigma(|x|)|x|\}.
                   \label{(i)}
\end{equation}
}
%

\section{Some Useful Tools}

\noindent \textbf{Monotonicity Formula}

As is common for these types of problems, the following monotonicity
formula will be very useful for us. For a function $v$ let us define
\begin{equation*}
I(r,v,x^0)=\int_{B(x^0,r)}\frac{|\nabla
v(x)|^2}{|x-x^0|^{n-2}}dx.
\end{equation*}

\begin{theorem}
\cite{1} Let $h_1,\ h_2$ be two non-negative continuous
sub-solutions of $\Delta u=0$ in $B(x^0,R)$. Assume further that
$h_1h_2=0$ and that $ h_1(x^0)=h_2(x^0)=0$. Then the following
function is monotone in $r\ (0<r<R)$
\begin{equation}
\phi(r,h_1,h_2,x^0)=\frac1{r^4}I(r, h_1,x^0)I(r, h_2,x^0).
                             \label{Fmono}
\end{equation}
Moreover, if any of the sets $supp(h_i)\cap\partial B(x^0,r)$
digresses from a spherical cap by a positive area, then either
$\phi^{\prime}(r)>0$ or $ \phi(r)=0$.
\end{theorem}

We use the abbreviated notations
$\phi(r)=\phi(r,h_1,h_2)=\phi(r,h_1,h_2,0)$. \vspace{2mm}

\noindent \textbf{Odd Reflection}\\
In order to be able to use the monotonicity formula, in some
cases we extend $P^+_r(M)$ to the class $P^{*}_r(M)$ of functions
that are defined in the whole $B_r$:
\begin{equation}
\tilde u(x_1,x_2,\ldots,x_n)=
\left\{
   \begin{array}{ll}
     u(x_1,x_2,\ldots,x_n) &  \text{for } x_1\geq 0,\\
     -u(-x_1,x_2,\ldots,x_n) &  \text{for } x_1<0,\\
    \end{array}
                    \right.   \label {remark-odd-reflection}
\end{equation}
We also set
$$\begin{array}{lll}
\Omega^-&= & \{(x_1,x_2,\ldots,x_n)\in B^-_r:
                         (-x_1,x_2,\ldots,x_n)\in\Omega\},\\
\Omega^{\prime}&= & \Omega\cup\Omega^-.
\end{array}
$$


\begin{definition}
A function $u$ (not identically zero) belongs to the class
$P^{*}_r(M)$, if $u$ satisfies:

\begin{enumerate}
\item $\Delta u=\chi_\Omega-\chi_{\Omega^-}\ \text{in}\ B_r$,
      in the sense of distributions ,

\item $|\nabla u|=0\ \text{in}\ B_r\backslash\Omega^{\prime}$,

\item $u=0\ \text{on}\ \Pi_r$,

\item $\|u\|_{\infty,B_r}\leq M$.
\end{enumerate}
\end{definition}

We also define $P^{*}_r(0,M)$ as a subclass of $P^{*}_r(M)$ for
which the origin belongs to the free boundary. 

\begin{lemma}
If $u$ is a solution of problem \eqref{pr}, then for all
$x^0\in\Gamma$ and $ 0<r<dist(x^0,\partial B)$ we have
\begin{equation}
\sup_{B^+(x^0,r)}u>u(x^0).
\end{equation}
\label{lemma-patch-changing}
\end{lemma}

\proof For $x^0\in\Gamma\backslash\Gamma^*$ we may apply the
strong maximum principle to $u$ in $\Omega\cap B^+(x^0,r)$ to
obtain the result ($u $ cannot be constant in $B^+(x^0,r)$ since
$x^0\in\partial\Omega$). Next let $x^0\in\Gamma^*$. We know that
$\Delta u\geq0$ in $B^+(x^0,r)$. If $ \sup_{B^+(x^0,r)}u\leq
u(x^0)$, then by Hopf's lemma type argument we have $
\frac{\partial u}{\partial x_1}(x^0)<~0$, which contradicts to
$|\nabla u|=0$ on $\Gamma^*$. Here we have used that $u\in
C^{1,\alpha}(\overline{B} ^+_1)$, for some $0<\alpha < 1$. \qed

\vspace{3mm}

\noindent \textbf{Blow-Up and Non-Degeneracy}\newline For a
function $u$, point $x^0\in\Gamma(u)$ and $r>0$ we consider
the following scaling
\begin{equation*}
u_r(x):=\frac{u(rx+x^0)-u(x^0)}{r^2}.
\end{equation*}

\begin{remark}
If $u\in P^+_r(x^0,M)$, $\|D_{ij}u\|\leq CM$ then
\begin{eqnarray*}
u_s(x)\in P^+_1(0,M),\\
\begin{array}{lll}
\Omega(u_s)&=&\Omega_s(u),\\
\Lambda(u_s)&=&\Lambda_s(u),\\
\Gamma(u_s)&=&\Gamma_s(u),
\end{array}
\end{eqnarray*}
where $E_s=\{x:sx+x^0\in E\}$ for any set $E$.
                      \label{remark-scalings}
\end{remark}
The uniform limit of $u_{r_j}$ when $r_j\to0$, is called
blow-up of $u$.

Obviously, if $u(x)=u(x^0)+o(|x-x^0|^2)$ then the blow-up
limit of $u$ will degenerate to be identically zero. The
following lemma shows, that in our problem we have
\textbf{non-degeneracy:}

\begin{lemma}
\cite{8} If $u\in P_1^+(0,M)$, $x^0\in\overline\Omega\cap
B_{1/2}$ such that $u(x^0)\geq0$, then
\begin{equation}
\sup_{B^+(x^0,r)}{u}\geq u(x^0)+Cr^2, \qquad \text{for all}
\ r<dist(x^0,\partial B_1),  \label{eq-non-deg}
\end{equation}
where $C=C(n)$. If $u(x^0)<0$ then \eqref{eq-non-deg} holds
with smaller $C_n$, provided $B(x^0,r)\subset B_1^+$.
                            \label{lemma-non-deg}
\end{lemma}
\proof We will consider the cases when $u(x^0)\geq0$ and
$u(x^0)<0$ separately. \\
\noindent {\sl Case 1.}
$u(x^0)\geq0$. We can assume that $x^0\in\Omega\cap
B_{1/2}$, since if \eqref{eq-non-deg} holds for all
$x^0\in\Omega\cap B_{1/2}$ then it will be true also for
all $x^0\in\overline{\Omega}\cap B_{1/2}$.

Let us set
\begin{equation}
v(x)=u(x)-u(x^0)-\frac{1}{2n}|x-x^0|^2 . \label{v(x)}
\end{equation}
There exists $x^1\in\overline{B^+}(x^0,r)$ such that the
following holds:
\begin{equation}
v(x^1)=\sup_{B^+(x^0,r)}{v}.  \label{v(x^1)}
\end{equation}
To prove this case, it is enough to prove the following two
steps:

\begin{itemize}
\item $v(x^1)\geq0$,

\item $x^1\in\partial B^+(x^0,r)\backslash\Pi(x^0,r)$.
\end{itemize}

\vspace{1mm} The first step simply follows from the fact
that $v(x^1)\geq v(x^0)=0$.\\
To prove the second step assume $x^1\in B^+(x^0,r)$. Then
from \eqref{v(x^1)} we have $|\nabla v|(x^1)=0$. Thus by
\eqref{v(x)}
\begin{equation}
(\nabla u)(x^1)=\frac{1}{n}(x^1-x^0). \label{nabla_u}
\end{equation}
Now, if $x^1\neq x^0$, then $(\nabla u)(x^1)\neq0$, i.e.,
$x_1\in\Omega$. But $\Delta v\geq0$ in $\Omega$ and \eqref{v(x^1)}
together with maximum principle gives us that $v(x)\equiv
constant=:C$ in $\overline{\Omega}\cap \overline{B^+(x^0,r)}$.
Particularly, $C=v(x^0)=0$ so we have
$u(x)=u(x^0)+\frac{1}{2n}|x-x^0|^2 $ and $(\nabla
u)(x)=\frac{1}{n}(x-x^0)$ in $\overline{\Omega}\cap
\overline{B^+(x^0,r)}$. Without loss of generality we may assume
that there exists $y\in\partial\Omega\cap {B^+(x^0,r)}$. Indeed,
if such an $y$ does not exist, then we have $B^+(x^0,r) \subset
\Omega$ and thus $u(x)=u(x^0)+\frac{1}{2n}|x-x^0|^2$ in
$B^+(x^0,r)$, which implies \eqref{eq-non-deg}.

 Thus we get $|\nabla u(y)|=\frac{1}{n}|y-x^0|\neq 0$, which is a contradiction,
 since $|\nabla u|=0$ on $\partial\Omega$.

If $x^1=x^0$, then again $x^1=x^0\in\Omega$ which contradicts to
\eqref{nabla_u}. Thus we have $x^1\in\partial B^+(x^0,r)$.

Finally, if $x^1\in\Pi(x^0,r)$, then because $u(x^0)\geq0$,
we get the following contradiction $0>v(x^1)\geq v(x^0)=0$.
\\
\noindent {\sl Case 2.}  $u(x^0)<0$. The proof of this
case is essentially the same as the proof of the previous
one, except we do not have
 $x^1\in\Pi(x^0,r)$, which was the only occasion  when we
used the nonnegativity of $u(x^0)$. \qed

\section{Proof of theorem A}

First let us extend $u$ from the class $P^+_1$ to the class
$P^{*}_1$ by odd reflection, as in
\eqref{remark-odd-reflection}. Set
\begin{equation*}
S_j(z,u)=\max_{B_{2^{-j}}(z)}|u(x)-u(z)| .
\end{equation*}
It is enough to prove the following lemma:

\begin{lemma}
There exist a constant $C_0$ depending only on $n$, such
that for every $u\in P^{*}_1(0,M)$, $j\in\mathbb{N}$ and
$z\in\Gamma(u)\cap B_{1/2}$
\begin{equation}
S_{j+1}(z,u)\leq max\{S_j(z,u)2^{-2},C_0M2^{-2j}\} .
\end{equation}
\label{LforN}
\end{lemma}

\proof

If the conclusion in the lemma fails, then there exist
sequences $\{u_j\}\subset~P^{*}_1(0,M),\
\{z_j\}\subset\Gamma(u_j)\cap B_{1/2},\
\{k_j\}\subset\mathbb{N},\ k_j\nearrow\infty$ such that
\begin{equation*}
   S_{k_{j}+1}(z_j,u_j)>max\{S_j(z_j,u_j)2^{-2},
   \ Mj2^{-2k_j}\}   \qquad \forall j\in\mathbb{N}.
\end{equation*}
Observe that $u_j\in P^{*}_1(0,M)$ implies
\begin{equation*}
\Delta u_j=\left\{
\begin{array}{ll}
\chi_{\Omega_j} & \text{if $x_1>0$}, \\
-\chi_{\Omega_j^-} & \text{if $x_1<0$}.%
\end{array}
\right.
\end{equation*}
Now consider the following scalings
\begin{equation*}
   \tilde{u}_j(x)=\frac{u_j(z_j+2^{-k_j}x)-u_j(z_j)}
   {S_{k_j+1}(z_j,u_j)} \qquad \text{in $B_1$}.
\end{equation*}
The following results can be obtained by computation like in
[14]:

\begin{itemize}
\item $\|\tilde{u}_j\|_{\infty,B}=\frac{S_{k_j}(z_j,u_j)}
      {S_{k_j+1}(z_j,u_j)} \leq4$ ,

\item $\|\tilde{u}_j\|_{\infty,B_{1/2}}=1$ ,

\item $\tilde{u}_j(0)=|\nabla\tilde{u}_j(0)|=0$,

\item $\|\Delta\tilde{u}_j\|_{\infty,B}\leq
      \frac4{j}\longrightarrow0,\ \text{when}\ j\to\infty$.
\end{itemize}

By compactness there exists a subsequence of
$\{\tilde{u}_j\}$ converging to a function $u_0$ in
$W^{2,p}(B_{1/2})\cap C^{1,\alpha}(B_{1/2})$. We have
$u_0(0)=|\nabla u_0(0)|=0$. For the renamed converging
subsequence $\tilde{u}_j$ set
\begin{equation*}
v=D_eu_0,\quad v_j=D_eu_j,\quad \tilde{v}_j=D_e\tilde{u}_j,
\end{equation*}
where $e$ is a fixed direction orthogonal to $e_1$. Obviously we
have that $v$ is the $C^{0,\alpha}$ limit of the sequence
$\tilde{v_j}$ in $B_1$, $v^{\pm}_j(0)=0$ and $\Delta v^{\pm}_j=0$.
Next we will use the monotonicity formula \eqref{Fmono} for the
sequence $\{v_j^\pm\}$ to get
\begin{equation}
     \frac1{r^{2n}}\int_{B_r}|\nabla v_j^+|^2
      \int_{B_r}|\nabla v_j^-|^2\leq C \quad \forall r,j,%
\end{equation}
where $C$ depends only on $M$. From here, using Poincare
inequality and letting $j$ go to infinity we obtain
\begin{equation}
\int_{B_1}|\tilde{v}^+-M^+|^2\int_{B_1}|\tilde{v}^--M^-|^2=0,%
                         \label{=0}
\end{equation}
where $M^\pm$ is the mean value of $v^\pm$ in $B_1$. Since
$v(0)=0$, from \eqref{=0} we have that either of $v^\pm$ is $0$.
Using the maximum principle we get $D_eu_0=v\equiv0$. This means
that $u_0$ depends only on the $x_1$ direction. Since it is
harmonic and has a second order growth, we obtain
$u_0(x)=ax_1+b$. Finally $u_0(0)=u_0^{\prime}(0)=0$ brings us to
the statement that $u_0\equiv0$, which contradicts
$\|\tilde{u}_j\|_{\infty,B_{1/2}}=1$, $\forall
j$. \qed \\

From this lemma we have the following inequality:
\begin{equation*}
|u(x)|\leq CMd(x)^2,\quad d(x)=dist(x,\partial\Omega),
\end{equation*}
for the points in $B_{1/2}$ that are close to $\Gamma$. This
together with elliptic estimates for the points close to
$\partial\Omega\cap\Pi$ gives us $\sup_{B^+_{1/2}}{|D_{ij}u|}\leq
CM$.
\begin{remark}
The free boundary has zero Lebesgue measure. \label{LebMes}
\end{remark}
This can be checked similarly as it is done in \cite{5}
 (see also [7; General remarks]). Only the non-degeneracy
 and $C^{1,1}$ properties of the solution are used in the
 proof.
\section{Proof of theorem B}

The first part of the proof consists of using the quadratic growth
of solutions to show that they are two dimensional. In the second
part we solve the problem in two dimensions.
\begin{lemma}\cite{14}
The global solutions are two dimensional.
\end{lemma}

Although the proof is very similar to what is found in \cite{14},
for readers convenience we will give an outline here. See
\cite{14} for more details.

At first we fix a direction $e$ orthogonal to $e_1$ and consider
$(D_eu)^\pm$. Since $(D_eu)^\pm$ vanish on $\Pi$, we can extend
them to the entire space $\mathbf{R}^n$ defining them as zero in
$\mathbf{R}^-_n$. Then, using the results of Theorem A, a
compactness argument and monotonicity formula \eqref{Fmono}, we
obtain that $D_eu$ doesn't change sign. Assume it is non-negative
(the non-positive case can be treated similarly), then by strong
maximum principle we get that on connected components of $\Omega$,
$D_eu$ must be strictly positive or identically zero. If $D_eu$ is
zero for all directions orthogonal to $e_1$, then $u$ is one
dimensional, so we have the representation b). If there is a
direction $e$ orthogonal to $e_1$ such that
\begin{equation}
D_eu>0,  \label{D_eu}
\end{equation}
then it can be proved that $u$ is two dimensional on every
connected component of $\Omega$. Thus it is enough to consider
the two dimensional problem. We treat two different cases.
\newline

\noindent \textbf{Case a)} When
$\overline{\Omega}=\overline{\mathbf{R}^n_+} $ . Then, since
$\partial\Omega$ has zero Lebesgue measure (see Remark
\ref{LebMes}), $D_2u$ is harmonic in upper half space and
vanishes on $\Pi$, so we can continue it harmonically by
reflection into entire space. Using Liouville's theorem and the
quadratic growth of $u$ we can conclude that $D_2u$ is linear.
Simple calculation then gives us the desired result.\newline

\noindent \textbf{Case b)} When $\overline{\Omega} \neq
\overline{\mathbf{R}^n_+} $. Then the interior of $\Lambda$ is
non-empty and we can take a ball $B(x^0,2R)\subset\Lambda(u)$.
Denote
\begin{equation*}
K(x^0,R):=\{(x_1,x_2-s):(x_1,x_2)\in B(x^0,R), s\geq0\}.
\end{equation*}
\noindent We claim
\begin{equation*}
\partial\Omega\cap K(x^0,R)=\emptyset.
\end{equation*}
Suppose this fails, and let $y\in\partial\Omega\cap K(x^0,R)$
(see Figure 3). Let also $u\equiv C_1$ in the connected component
of $\Lambda(u)$ that contains $x^0$. From $D_2u>0$ (see
\eqref{D_eu}) follows that $u\leq C_1$ in $K(x^0,2R)$. Using the
strong maximum principle ($u$ is subharmonic) and that $u(y)=C_1$
we conclude that $u=C_1$ in $B(y,R)$, which contradicts the
assumption $y\in\partial\Omega$. Hence
$K(x^0,R)\subset~\Lambda(u)$.

In order to prove that $u$ is one dimensional, let us
extend $u$ to a function $\tilde u$ defined in the whole
space $\mathbf{R}^2$ as in \eqref{remark-odd-reflection}.
 Then for $D_2\tilde{u}$ we will have:
\begin{equation}
\begin{array}{ll}
D_2\tilde{u}(x_1,x_2)=
\left\{
\begin{array}{ll}
D_2u(x_1,x_2) & \text{in}\ \mathbf{R}_+^2, \\
-D_2u(-x_1,x_2) & \text{in}\ \mathbf{R}_-^2. \\
\end{array}
\right.%
\end{array}
\label{eq-D2u-odd}
\end{equation}
We next consider the blow-up of $\tilde{u}$ at $\infty$:
$\tilde{u}_\infty(x)=\lim_{r_j\to\infty}\tilde{u}_{r_j}(x)$
where, as usual, $r_j \nearrow\infty$ and
$\tilde{u}_{r_j}(x)={\tilde{u} (r_jx)}/{r^2_j}$. Also
observe that by the definition of global solutions
$\tilde{u}_{r_j}$ is bounded. Writing the monotonicity
formula for functions $(D_2\tilde{u})^\pm$  we get:
\begin{equation}
\phi(r,D_2\tilde{u})\leq\phi(r_j,D_2\tilde{u})\leq
\lim_{r_j\to\infty}\phi(r_j,D_2\tilde{u})=
\phi(1,D_2\tilde{u}_\infty)=C
\label{eq-fi-chain}
\end{equation}
for $0<r<r_j$. Next, let us observe that for a fixed $s>0$
the following holds
\begin{equation*}
C=\lim_{r_j\to\infty}\phi(sr_j,D_2\tilde{u})=
\lim_{r_j\to\infty}\phi(s,D_2\tilde{u}_{r_j})=
\phi(s,D_2\tilde{u}_\infty).
\end{equation*}
 Hence $\phi(s,D_2\tilde{u}_\infty)=C=constant$ for all
 $s>0$.

In order to complete the proof of the theorem we need the
following lemma.
\begin{lemma}
For any $s>0$ we have
\begin{equation*}
\phi(s,D_2\tilde{u}_\infty)=C=0.
\end{equation*}
       \label{s>0}
\end{lemma}
\textit{Proof.} We prove the lemma by a contradictory argument.
Let us assume there exists $s>0$ such that
$\phi(s,D_2\tilde{u}_\infty)=C$ and $C\neq0$. First we observe
that $D_2\tilde{u}_\infty>0$ in $\mathbf{R}^2_+$. Since otherwise
(by maximum principle) there exists a ball
$B(y_0,t)\subset\{\mathbf{R} ^2_+\setminus
supp(D_2\tilde{u}_\infty)^+\}$, $t>0$. Also we have that $
supp(D_2\tilde{u}_\infty)^+\cap\partial B_{|y_0|}(0)$ digresses
from a spherical cap by a positive area. Then the monotonicity
formula applied for $(D_2\tilde{u})^\pm$ on the ball
$B_{|y_0|}(0)$ implies that $\phi(|y_0|,D_2\tilde{u}
_\infty)=C=0$, which is a contradiction.
\begin{figure}[tbp]
\begin{picture}(300,150)(0,0)
\put(10,-3){\includegraphics[height=5cm]{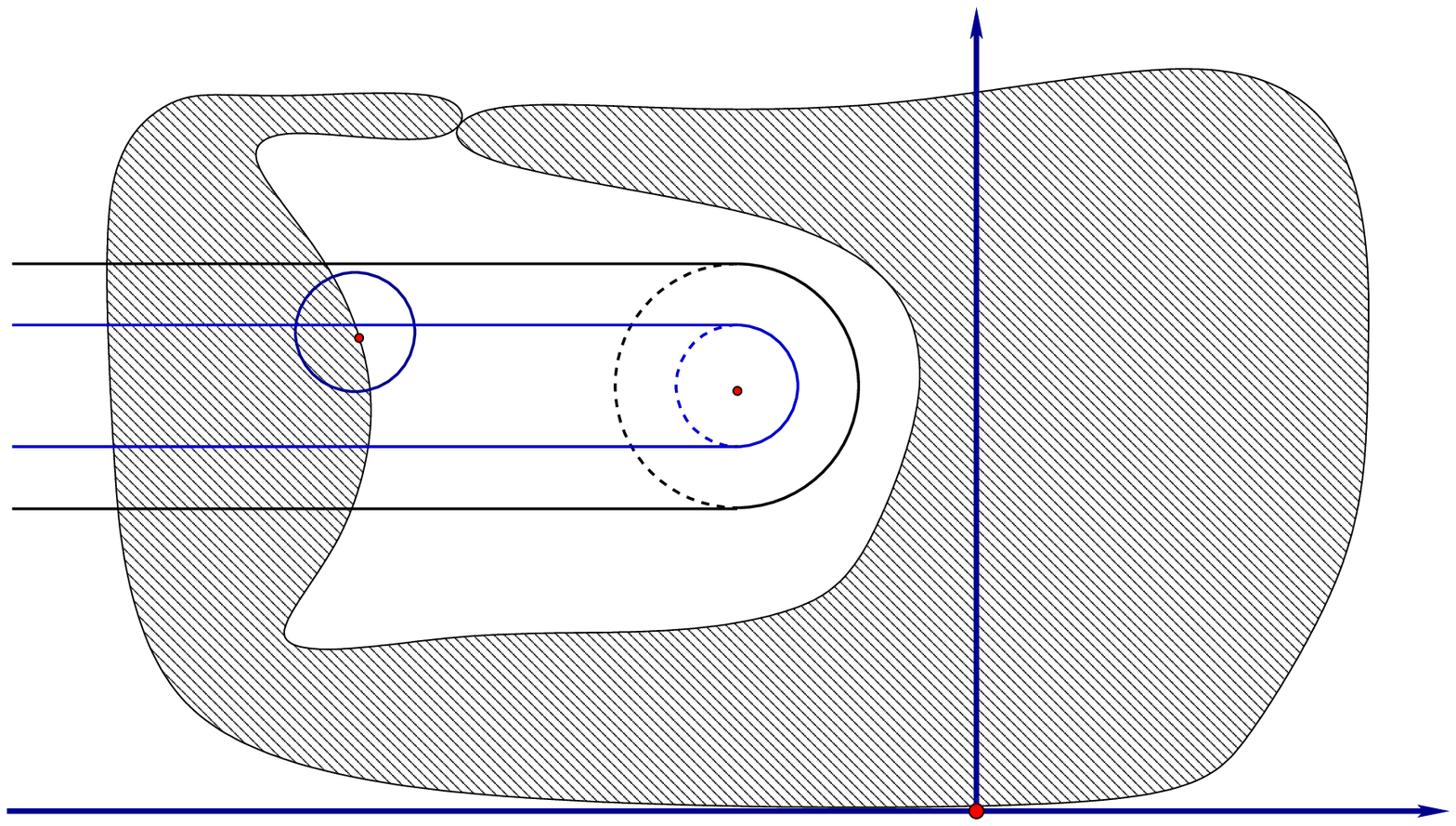}}%
\put(150,15){\Large$\Omega$}%
\put(75,54){\tiny$K(x^0,2R)$}%
\put(182,135){$x_1$}%
\put(75,38){$\Lambda(u)$}%
\put(245,-14){$\Pi$}%
\put(73,77){\tiny$y$}%
\put(138,68){\tiny$x^0$}%
\end{picture}
\caption{}
\end{figure}
Also notice that $suppD_2\tilde{u}_\infty=\overline{\mathbf{
R}^2_+}$ implies $\Omega(\tilde{u}_\infty)=\mathbf{R}^2_+$. Just
like in case a), we have that
$\tilde{u}_\infty(x)=\frac{x_1^2}2+ax_1x_2+\alpha x_1$, where
$a,\alpha\in\mathbf{R}$. Also $D_2\tilde{u}_\infty>0$ implies that
$ a>0 $. To get a contradiction, it is enough to prove that
$D_1\tilde{u} _\infty$ is zero at two different points
$(0,x^{\prime})$, $ (0,x^{\prime\prime})$ of $\Pi$. Indeed, assume
$D_1\tilde{u} _\infty(0,x_2^{\prime})=0$ and
$D_1\tilde{u}_\infty(0,x_2^{\prime\prime})=0$, then we have
$ax^{\prime}=\alpha=ax^{\prime\prime}$, and the only possibility
is that $a=0$, contradicting $D_2\tilde{u}>0$.

Let us now show that there exists two different points on $\Pi$,
where $D_1\tilde{u}_\infty$ is equal to zero. In fact one can
prove even more. Namely $|\nabla \tilde{u}|$ vanishes on
$\Pi\cap\{x_2<0\}$. Recall that we have $B(x^0,R)\subset
K(x^0,R)\subset\Lambda(\tilde{u})$. \noindent Denote
\begin{equation*}
K_{r_j}:=K\left(\frac{1}{r_j}x^0,\frac{R}{r_j}\right)\quad
\text{and} \quad l_j:=\left\{\frac{1}{r_j}x^0-se_2,\ e_2=(0,1),\
s>0\right\}.
\end{equation*}
Fix an $s>0$ and consider the sequence $y_j:=x^0/r_j-se_2$.
Obviously, $ y_j\in l_j\subset
K_{r_j}\subset\Lambda(\tilde{u}_{r_j})$.
 Recalling that $\tilde{u}_{r_j}$ converges to $u_\infty$ in
$(W^{2,p}_{loc}\cap
C^{1,\alpha}_{loc})(\mathbf{R}_+^n\cup\Pi)$,
 we get
\begin{equation*}
|\nabla \tilde{u}_\infty(0,-s)|=0, \ \forall s>0.
\end{equation*}
This completes the proof of the lemma. \qed
\\

Finally, using \eqref{eq-fi-chain}, \eqref {eq-D2u-odd} and
positivity of $D_2u$ in $R_+^2$ we have
\begin{eqnarray}
0=\phi(1,D_2\tilde{u}_\infty)\geq\phi(r,D_2\tilde{u})=
   \frac1{r^4}I\left(r,(D_2\tilde{u})^+,x^0\right)I
   \left(r,(D_2\tilde{u})^-,x^0\right)\notag \\
=\frac1{r^4}\left(\int_{B_r^+}
 \frac{|\nabla(D_2u)(x_1,x_2)|^2dx}{|x|^{n-2}}\right)
 \left(\int_{B_r^-}\frac{|\nabla(D_2u)(-x_1,x_2)|^2dx}
 {|x|^{n-2}}\right)  \notag \\
= \frac1{r^4}\left(\int_{B_r^+}\frac{|\nabla(D_2u)|^2dx}
{|x|^{n-2}}\right)^2,\notag
\end{eqnarray}
for any $r>0$. This gives us $|\nabla D_2\tilde{u}| \equiv0$ in
$\mathbf{R}^2$. Hence $D_2u\equiv constant=D_2u(0)=0$. Therefore
we get $u$ is one dimensional. Simple calculations combined with
$ C^{1,\alpha}$ regularity of the solutions accomplish the proof
of the theorem. \qed

\section{Proof of theorem C}

It is enough to check that for every given $\epsilon$ there exists
$\rho=\rho_\epsilon$ such that for all
$x^0\in\partial\Omega\cap B^+_{\rho_\epsilon}$%
\begin{equation}
x^0\in B^+_{\rho_\epsilon}\backslash K_\epsilon,  \label{tang}
\end{equation}
where $K_\epsilon=\{x:x_1>\epsilon(x_2^2+\ldots+x_n^2)^{1/2}\}$.
Then we may choose $r_0=\rho_{\{\epsilon=1\}}$ and $\sigma$ given
by the inverse of $ \epsilon\to\rho_\epsilon$.

 Conversely, suppose that \eqref{tang}
fails, then there exists a sequence $ u_j\in P_1^+(0,M)$,
$x^j\in\partial\Omega(u_j)\cap B^+_{\rho_j}$ such that
$\rho_j\to0$ and $x^j\in B^+_{\rho_j}\cap\overline{K_\epsilon}$.
Now, for every scaled function
$\tilde{u}_j(x)=u_j(x|x^j|)/|x^j|^2$ we have a point $
\tilde{x}^j\in K_\epsilon$. There exists converging subsequences
of $\tilde{u}_j\to u_0$ and $\tilde{x}^j\to x^0$ such that
$x^0\in\overline{K}_\epsilon\cap\partial B_1$. Since
$x^0\in\Gamma$, $0\in\Gamma$, and $u_0$ is a global solution we
have a contradiction to Theorem B. \qed
\newline

\section*{Acknowledgements}
I would like to express a deep gratitude to my Ph.D. thesis advisor
Henrik Shahgholian for suggesting the problem and for valuable
discussions.


\begin{thebibliography}{1}
\bibitem[1]{1} Alt, H. W.; Caffarelli, L. A.; Friedman, A.
 Variational problems with two phases and their free
boundaries. Trans. AMS {\bf 1984}, \textit{282}(4-5), 431-461.

\bibitem[2]{2}  Athanasopoulos, I.; Caffarelli, L. A.; Salsa S.
Regularity of the free boundary in parabolic phase-transition
problems. Acta Math. {\bf 1996}, \textit{176}, 245-282.

\bibitem[3]{3} Caffarelli, L. A. The regularity of free
boundaries in higher dimension. Acta Math. {\bf 1977},
\textit{139}, 155-184.

\bibitem[4]{4} Caffarelli, L. A. Compactness methods in free
boundary problems. Comm. P.D.E. {\bf 1980}, \textit{5}, 427-448.

\bibitem[5]{5} Caffarelli, L. A. The obstacle problem revisited.
J.~Fourier Anal. Appl. {\bf 1998}, \textit{4}(4-5), 383-402.

\bibitem[6]{6} Chapman, S. J. A mean-field model of superconducting
vortices in three dimensions. SIAM J. App. Math. {\bf 1995},
\textit{55}, 1259-1274.

\bibitem[7]{7} Caffarelli, L. A.; Karp, L.; Shahgholian, H.
Regularity of a free boundary with application to the Pompeiu
problem. Ann. Math. {\bf 2000}, \textit{151}, 269-292.

\bibitem[8]{8} Caffarelli, L. A.; Salazar, J.Solutions of fully
nonlinear elliptic equations with patches of zero gradient:
existence, regularity and convexity of level curves. Trans. Amer.
Math. Soc. {\bf 2002}, \textit{354}(8),  3095-3115 (electronic).

\bibitem[9]{9} Caffarelli, L. A.; Salazar, J.; Shahgholian, H.
Free boundary regularity for a problem arising in
Superconductivity.  Arch. Ration. Mech. Anal. {\bf 2004}, \textit
{171}(1), 115-128.

\bibitem[10]{10} Elliott, C. M.; Sch{\"a}tzle, R.; Stoth, B.E.E.
Viscosity solutions of a degenerate parabolic-elliptic system
arising in the mean-field theory, Arch. Rat. Mech. Anal. {\bf
1998}, \textit{145}, 99-127.

\bibitem[11]{11} Friedman, A. \textit{Variational principles and
free-boundary problems}, 2nd Ed.; Robert E. Krieger Publishing
Co., Inc.: Malabar, FL, 1988.

\bibitem[12]{12} Gilbarg, D.; Trudinger, N. S.; \textit{Elliptic Partial
Differential Equations of Second Order}, 2-nd Ed.;
Springer-Verlag: Berlin Heidelberg New York Tokyo, 1983.

\bibitem[13]{13} Karp, L.; Shahgholian, H. Regularity of a free
boundary problem. J. Geom. Anal. {\bf 1999}, \textit{9}(4),
653-669.

\bibitem[14]{14} Shahgholian, H.; Uraltseva, N. N. Regularity
properties of a free boundary near contact points with the fix
boundary. Duke math. J. {\bf 2003}, \textit{116}(1), 1-34 .

\bibitem[15]{15} Uraltseva, N. N. $C^1$-regularity of the boundary of
the noncoincidence set in the obstacle problem, (in Russian).
Algebra \& Analysis {\bf 2003}, \textit{8}(2). English
translation St. Petersburg Math. J. {\bf 1997}, \textit{8},
341-353.

\bibitem[16]{16} N. N. Uraltseva, \textit{On the properties of a free
boundary on a neighborhood of the points of contact with the
known boundary} , (in Russian). Zap. Nauchn. Sem. S.-Peterburg.
Otdel. Mat. Inst. Steklov. (POMI) {\bf 1997} \textit{249}.
English translation in  J. Math. Sci. (New York) {\bf 2000},
\textit{101}(5), 3570-3576.



\end{thebibliography}
\end{document}